\documentclass[12pt]{article}
\usepackage[letterpaper, left=1.5in, right=1in,top=1in,bottom=1in]{geometry}
\usepackage{tikz}
\usepackage{chngcntr}
\counterwithin{figure}{section}
\usetikzlibrary{arrows}
\usepackage{placeins}	
\usepackage{amsfonts}
\usepackage{graphicx}
\usepackage{standalone}
\usepackage{amsmath, amsthm, amssymb}
\usepackage[capitalize, noabbrev]{cleveref}
\newtheorem{thrm}{Theorem}[section]

\newtheorem{lem}[thrm]{Lemma}
\newtheorem{cor}[thrm]{Corollary}

\newtheorem{ques}[thrm]{Question}

\newenvironment{pf}           {\noindent{\bf Proof:} }%
                                {\null\hfill$\Box$\par\medskip}
                           
\usepackage[nottoc]{tocbibind}

\begin{document}

\title{Orthogonal Colourings of Cayley Graphs}
\author{Jeannette Janssen\thanks{Email: Jeannette.janssen@dal.ca} , Kyle MacKeigan\thanks{Email: Kyle.m.mackeigan@gmail.com}}
\maketitle

\begin{abstract}
Two colourings of a graph are orthogonal if they have the property that when two vertices are coloured with the same colour in one colouring, then those vertices receive distinct colours in the other colouring. In this paper, orthogonal colourings of Cayley graphs are discussed. Firstly, the orthogonal chromatic number of cycle graphs are completely determined. Secondly, the orthogonal chromatic number of certain circulant graphs is explored. Lastly, orthogonal colourings of product graphs and Hamming graphs are studied.
\end{abstract}

\section{Introduction}
Two colourings of a graph are \textit{orthogonal} if they have the property that when two vertices are coloured with the same colour in one colouring, then those vertices must have distinct colours in the other colouring. A $k$-\textit{orthogonal colouring} is a collection of $k$ mutually orthogonal colourings. For simplicity, a $2$-orthogonal colouring is called an \textit{orthogonal  colouring}. The \textit{$k$-orthogonal chromatic number} of a graph $G$, denoted $O\chi_k(G)$, is the minimum number of colours required for a $k$-orthogonal colouring. The $2$-orthogonal chromatic number is simply called the \textit{orthogonal chromatic number}, and is denoted $O\chi(G)$.  

Orthogonal colourings were first defined in 1985 by Archdeacon, Dinitz, and Harary, under the context of edge colourings \cite{archdeacon1985orthogonal}. Later in 1999, Caro and Yuster revisited orthogonal colourings, this time under the context of vertex colourings \cite{caro1999orthogonal}. Then in 2013, Ballif studied upper bounds on sets of orthogonal vertex colourings \cite{ballif2013upper}.

Let $(G,\circ)$ denote a group $G$ with group operation $\circ$ and let $S$ be a generating set of $G$. Then, the associated \textit{Cayley graph}, denoted $\Gamma(G,S)$, has a vertex for each element of $G$ and there is a directed edge between two elements $u$ and $v$ if and only if $u\circ v^{-1}\in S$. It is assumed that $S$ is closed under inverses and that the additive identity is not in $S$. This is so that the Cayley graphs considered in this paper are simple, and undirected. 

To start, the orthogonal chromatic number of cycle graphs are determined. Then, orthogonal colourings of more general circulant graphs are studied. Lastly, orthogonal colourings of Cartesian product graphs and Hamming graphs are explored.

\section{Orthogonal Colourings of Cycle Graphs}

The Cayley graphs of cyclic groups are called circulant graphs. For example, consider $\mathbb{Z}_n$, the integers modulo $n$ with group operation addition modulo $n$. Then, the Cayley graph $\Gamma(\mathbb{Z}_n,\{1,-1\})\cong C_n$, the cycle graph on $n$ vertices. For a graph $G$ with $n$ vertices, a lower bound for $O\chi(G)$ is $\lceil\sqrt{n}\,\rceil$. In this section, it is shown that $O\chi(C_n)=\lceil\sqrt{n}\,\rceil$ if and only if $n>4$. The following lemma establishes that this is the correct orthogonal chromatic number for most cases.

\begin{lem}\label{Lemma:Cycle Part 1}
For $n>4$, if $\lceil\sqrt{n}\,\rceil\nmid (n-1)$, and $\lceil\sqrt{n}\,\rceil\nmid \Big( n-1+\Big\lfloor\frac{n-1}{\lceil\sqrt{n}\,\rceil}\Big\rfloor\Big)$, then $O\chi(C_n)=\lceil\sqrt{n}\,\rceil$.
\end{lem}
\begin{pf}
Let $v_i\in \mathbb{Z}_n$ where $0\leq i <n$ and let $N=\lceil\sqrt{n}\,\rceil$. Set $c_1(v_i)=i(\textrm{mod}~N)$ and $c_2(v_i)=\big(i+\big\lfloor\frac{i}{N}\big\rfloor\big)(\textrm{mod}~N)$. For illustration, these two colourings are applied to $C_9$ in Figure \ref{Figure: Orthogonal colouring of C9}. Displayed next to each vertex $v_i$ is the colour pair $(c_1(v_i),c_2(v_i))$.
\begin{figure}[h!]
\centering
\begin{tikzpicture}[scale=.6,line cap=round,line join=round,>=triangle 45,x=1.0cm,y=1.0cm]
\draw [line width=2.pt] (0.52,2.82)-- (1.22,4.18);
\draw [line width=2.pt] (1.22,4.18)-- (2.96,4.66);
\draw [line width=2.pt] (2.96,4.66)-- (4.7,4.18);
\draw [line width=2.pt] (4.7,4.18)-- (5.5,2.7);
\draw [line width=2.pt] (5.5,2.7)-- (5.16,1.26);
\draw [line width=2.pt] (5.16,1.26)-- (4.,0.);
\draw [line width=2.pt] (4.,0.)-- (2.,0.);
\draw [line width=2.pt] (2.,0.)-- (0.8,1.38);
\draw [line width=2.pt] (0.8,1.38)-- (0.52,2.82);
\draw (2.5,4.5) node[anchor=north west] {$v_0$};
\draw (3.9,4.15) node[anchor=north west] {$v_1$};
\draw (4.4,3.2) node[anchor=north west] {$v_2$};
\draw (4.2,2) node[anchor=north west] {$v_3$};
\draw (3.2,0.9) node[anchor=north west] {$v_4$};
\draw (1.7,0.9) node[anchor=north west] {$v_5$};
\draw (.7,2) node[anchor=north west] {$v_6$};
\draw (0.55,3.2) node[anchor=north west] {$v_7$};
\draw (1.1,4.2) node[anchor=north west] {$v_8$};
\draw (2,5.8) node[anchor=north west] {(0,0)};
\draw (4.1,5.35) node[anchor=north west] {(1,1)};
\draw (5.5,3.25) node[anchor=north west] {(2,2)};
\draw (5,1.4) node[anchor=north west] {(0,1)};
\draw (3.5,0) node[anchor=north west] {(1,2)};
\draw (.9,0) node[anchor=north west] {(2,0)};
\draw (-.8,1.4) node[anchor=north west] {(0,2)};
\draw (-1.3,3.35) node[anchor=north west] {(1,0)};
\draw (-.1,5.35) node[anchor=north west] {(2,1)};
\begin{scriptsize}
\draw [fill=black] (2.,0.) circle (2.5pt);
\draw [fill=black] (4.,0.) circle (2.5pt);
\draw [fill=black] (2.96,4.66) circle (2.5pt);
\draw [fill=black] (0.8,1.38) circle (2.5pt);
\draw [fill=black] (5.16,1.26) circle (2.5pt);
\draw [fill=black] (0.52,2.82) circle (2.5pt);
\draw [fill=black] (5.5,2.7) circle (2.5pt);
\draw [fill=black] (1.22,4.18) circle (2.5pt);
\draw [fill=black] (4.7,4.18) circle (2.5pt);
\end{scriptsize}
\end{tikzpicture}
\caption{Orthogonal Colouring of $C_9$}
\label{Figure: Orthogonal colouring of C9}
\end{figure}
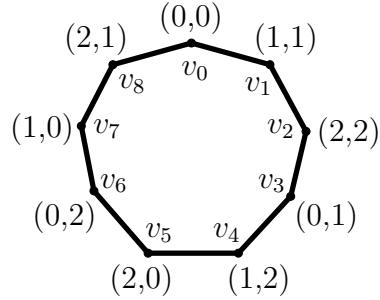

It is now shown that $c_1$ and $c_2$ are both proper colourings. For $0\leq i \leq n-2$, $c_1(v_i)=i(\textrm{mod}~N)\not\equiv (i+1)(\textrm{mod}~N)=c_1(v_{i+1})$. Then, since $N\nmid (n-1)$ by assumption, $c_1(v_{n-1})=(n-1)(\textrm{mod}~N)\neq 0 =c_1(v_0)$. Thus, $c_1$ is a proper colouring. To show that $c_2$ is a proper colouring, notice that $\big\lfloor\frac{i}{N}\big\rfloor\leq\big\lfloor\frac{i+1}{N}\big\rfloor\leq \big\lfloor\frac{i}{N}\big\rfloor+1 $. Therefore, for $0\leq i \leq n-2$, it follows that $1\leq c_2(v_{i+1})-c_2(v_i)\leq 2$. Thus, since $N>2$, $c_2(v_i)\neq c_2(v_{i+1})$. Now, since $N\nmid \big(n-1+\big\lfloor\frac{n-1}{N}\big\rfloor\big)$  by assumption, it follows that $c_2(v_{n-1})=\big(n-1+\big\lfloor\frac{n-1}{N}\big\rfloor\big)(\textrm{mod}~N)\neq 0 =c_2(v_0)$. Therefore, $c_2$ is a proper colouring. 

It remains to show that $c_1$ and $c_2$ are orthogonal colourings. Suppose otherwise, that is, $c_1(v_i)=c_1(v_j)$ and $c_2(v_i)=c_2(v_j)$ where $i\neq j$. Since $c_1(v_i)=c_1(v_j)$, this implies that $i=j+mN$ where $0<m<N$. Therefore:
\begin{align*}
c_2(v_i)&=\Big(j+\Big\lfloor\frac{j+mN}{N}\Big\rfloor\Big)(\textrm{mod}~N)\\
&=\Big(j+m+\Big\lfloor\frac{j}{N}\Big\rfloor\Big)(\textrm{mod}~N)\\
&=(c_2(v_j)+m) (\textrm{mod}~N).
\end{align*}

Since $c_2(v_i)=c_2(v_j)$, this gives that $m\equiv 0(\textrm{mod}~N)$, contradicting $0<m<N$. Therefore, $c_1$ and $c_2$ are orthogonal colourings of $C_n$. Since $c_1$ and $c_2$ both used $N$ colours, the fewest possible, $O\chi(C_n)=N$.
\end{pf}

Notice that the orthogonality property of $c_1$ and $c_2$ in Lemma \ref{Lemma:Cycle Part 1} did not depend on the assumed divisibility conditions. Therefore, the problem with using $c_1$ and $c_2$ in the cases where $N\mid (n-1)$ and $N\mid \big(n-1+\big\lfloor\frac{n-1}{N}\big\rfloor\big)$, is that there is some sort of colour conflict. The following lemma shows that this conflict can be resolved by assigning $v_{n-1}$ a different colour pair.

\begin{lem}\label{Lemma: Cycle Part 2}
For $n> 16$, if $\lceil\sqrt{n}\,\rceil\mid (n-1)$ or $\lceil\sqrt{n}\,\rceil\mid \Big(n-1+\Big\lfloor\frac{n-1}{\lceil\sqrt{n}\,\rceil}\Big\rfloor\Big)$, then $O\chi(C_n)=\lceil\sqrt{n}\,\rceil$.
\end{lem}
\begin{pf}
Let $N=\lceil\sqrt{n}\,\rceil$, $c_1(v_i)=i(\textrm{mod}~N)$ and $c_2(v_i)=\big(i+\big\lfloor\frac{i}{N}\big\rfloor\big)(\textrm{mod}~N)$. There are four cases to consider: Case 1: $N\mid (n-1)$ and $N\nmid\big(n+\big\lfloor\frac{n}{N}\big\rfloor\big)$. Case 2: $N\nmid n$ and $N\mid\big(n-1+\big\lfloor\frac{n-1}{N}\big\rfloor\big)$. Case 3: $N\mid (n-1)$ and $N\mid\big(n+\big\lfloor\frac{n}{N}\big\rfloor\big)$. Case 4: $N\mid\big(n-1+\big\lfloor\frac{n-1}{N}\big\rfloor\big)$ and $N\mid n$. For Case 1 and Case 2, the orthogonal colouring $(\bar c_1, \bar c_2)$ is used, where $\bar c_1$ and $\bar c_2$ are defined as follows:
\begin{align*}
\bar c_1(v_i)&= \begin{cases} 
      c_1(v_i) & 0\leq i\leq n-2 \\
      n(\textrm{mod}~N) & i=n-1 \\
   \end{cases}\quad\\
\bar c_2(v_i)&= \begin{cases} 
       c_2(v_i) & 0\leq i\leq n-2 \\
       \big(n+\big\lfloor\frac{n}{N}\big\rfloor\big)(\textrm{mod}~N) & i=n-1 \\
   \end{cases}
\end{align*}

It is shown that $\bar c_1$ and $\bar c_2$ are proper colourings in both cases. For $0\leq i\leq n-2$, $\bar c_1(v_i)=c_1(v_i)$ and $\bar c_2(v_i)=c_2(v_i)$. Therefore, by the proof of Lemma \ref{Lemma:Cycle Part 1}, there are no colour conflicts between these vertices. Note that in both Case 1 and Case 2, $n\not\equiv 0(\textrm{mod}~N)$ and $\big(n+\big\lfloor\frac{n}{N}\big\rfloor\big)(\textrm{mod}~N)\not\equiv 0(\textrm{mod}~N)$. Therefore, $\bar c_1(v_{n-1})\neq 0=\bar c_1(v_0)$, and $\bar c_2(v_{n-1})\neq 0=\bar c_2(v_0)$. 

Notice that $\big\lfloor\frac{n-2}{N}\big\rfloor\leq \big\lfloor\frac{n}{N}\big\rfloor\leq \big\lfloor\frac{n-2}{N}\big\rfloor+1$. Therefore, $2\leq \bar c_2(v_{n-1})-\bar c_2(v_{n-2})\leq 3$. Since $N>4$ by assumption, this implies that $\bar c_2(v_{n-1})\neq \bar c_2(v_{n-2})$. Also, since $N>4$, $\bar c_1(v_{n-2})=(n-2)(\textrm{mod}~N)\neq n(\textrm{mod}~N)=\bar c_1(v_{n-1})$. Therefore, $\bar c_1$ and $\bar c_2$ are proper colourings of $C_n$. It remains to show that $\bar c_1$ and $\bar c_2$ are orthogonal colourings.

 For $0\leq i \leq n-2$, there are no orthogonal conflicts on the vertices $v_i$ by the proof of Lemma \ref{Lemma:Cycle Part 1}. In Case 1, since $n\equiv 1(\textrm{mod}~N)$, the colour pair $(1,(1+\lfloor\frac{n}{N}\rfloor)(\textrm{mod}~N))$ is assigned to $v_{n-1}$. Let $i\equiv 1 (\textrm{mod}~N)$ and $i<n$. Let $m_1$ and $m_2$ be integers so that $i=m_1N+1$ and $n=m_2N+1$. Since $i<n$, this gives that $m_1<m_2$. Therefore, $\bar c_2(v_i)=(1+\lfloor\frac{m_1N+1}{N}\rfloor)(\textrm{mod}~N))=(1+m_1)(\textrm{mod}~N)\neq (1+m_2)(\textrm{mod}~N)=\bar c_2(v_{n-1})$. 
 
In Case 2, since $N\mid\big(n-1+\big\lfloor\frac{n-1}{N}\big\rfloor\big)$, the colour pair $(n(\textrm{mod}~N),1)$ is assigned to $v_{n-1}$. A similar argument as in Case 1 shows that there are no orthogonal conflicts. Therefore, in both Case 1 and Case 2, $\bar c_1$ and $\bar c_2$ are proper orthogonal colourings. For Case 3 and Case 4, a different orthogonal colouring $(\hat c_1,\hat c_2)$ is used, where $\hat c_1$ and $\hat c_2$ are defined as follows:
\begin{align*}
\hat c_1(v_i)&= \begin{cases} 
      c_1(v_i) & 0\leq i\leq n-2 \\
      (n+1)(\textrm{mod}~N) & i=n-1 \\
   \end{cases}\quad\\
\hat c_2(v_i)&= \begin{cases} 
       c_2(v_i) & 0\leq i\leq n-2 \\
       \big(n+1+\big\lfloor\frac{n+1}{N}\big\rfloor\big)(\textrm{mod}~N) & i=n-1 \\
   \end{cases}
\end{align*}

It is shown that $\hat c_1$ and $\hat c_2$ are proper colourings in both cases. For $0\leq i \leq n-2$, $\hat c_1(v_i)=c_1(v_i)$ and $\hat c_2(v_i)=c_2(v_i)$. Therefore, by the proof of Lemma \ref{Lemma:Cycle Part 1}, there are no colour conflicts between these vertices. Note that in both Case 3 and Case 4, $(n+1)\not\equiv 0(\textrm{mod}~N)$ and $\big(n+1+\big\lfloor\frac{n+1}{N}\big\rfloor\big)\not\equiv 0(\textrm{mod}~N)$ . Thus, $\hat c_1(v_{n-1})\neq 0=\hat c_1(v_0)$, and $\hat c_2(v_{n-1})\neq 0 =\hat c_2(v_0)$. 

Notice that $\big\lfloor\frac{n-2}{N}\big\rfloor\leq \big\lfloor\frac{n+1}{N}\big\rfloor\leq \big\lfloor\frac{n-2}{N}\big\rfloor+1$. Therefore, $3\leq \hat c_2(v_{n-1})-\hat c_2(v_{n-2})\leq 4$. Since $N>4$ by assumption, this implies that $\hat c_2(v_{n-2})\neq \hat c_2(v_{n-1})$. Also, since $N>4$, $\hat c_1(v_{n-2})=(n-2)(\textrm{mod}~N)\not\equiv (n+1)(\textrm{mod}~N)=\hat c_1(v_{n-1})$. Therefore, $\hat c_1$ and $\hat c_2$ are proper colourings of $C_n$. It remains to show that $\hat c_1$ and $\hat c_2$ are orthogonal colourings.

 For $0\leq i \leq n-2$, there are no orthogonal conflicts on the vertices $v_i$ by the proof of Lemma \ref{Lemma:Cycle Part 1}. In Case 3, since $n\equiv 1(\textrm{mod}~N)$, the colour pair $(2,(2+\lfloor\frac{n+2}{N}\rfloor)(\textrm{mod}~N))$ is assigned to $v_{n-1}$. Let $i\equiv 2(\textrm{mod}~N)$ and $i<n$. Let $m_1$ and $m_2$ be integers so that $i=m_1N+2$ and $n=m_2N+1$. Since $i<n$, this gives that $m_1<m_2$. Therefore, $\hat c_2(v_i)=(2+\lfloor\frac{m_1N+1}{N}\rfloor)(\textrm{mod}~N)=(2+m_1)(\textrm{mod}~N)\neq (2+m_2)(\textrm{mod}~N)=\hat c_2(v_{n-1})$. 
 
In Case 4, since $N\mid\big(n-1+\big\lfloor\frac{n-1}{N}\big\rfloor\big)$, the colour pair $((n+2)(\textrm{mod}~N),2)$ is assigned to $v_{n-1}$. A similar argument as in Case 3 shows that there are no orthogonal conflicts. Therefore, in both Case 3 and Case 4, $\hat c_1$ and $\hat c_2$ are proper orthogonal colourings.
\end{pf}

\begin{figure}[h!]
\centering
\begin{tikzpicture}[scale=.75,line cap=round,line join=round,>=triangle 45,x=1.0cm,y=1.0cm]
\draw [line width=2.pt] (-1,3.)-- (1,3.94);
\draw [line width=2.pt] (1,3.94)-- (3.,3.);
\draw [line width=2.pt] (3.,3.)-- (3.,1.);
\draw [line width=2.pt] (3.,1.)-- (1,0.);
\draw [line width=2.pt] (1,0.)-- (-1,1.);
\draw [line width=2.pt] (-1,1.)-- (-1,3.);
\draw (.3,4.85) node[anchor=north west] {(0,0)};
\draw (.3,-.1) node[anchor=north west] {(0,1)};
\draw (1.6,3.15) node[anchor=north west] {(1,1)};
\draw (1.6,1.7) node[anchor=north west] {(2,2)};
\draw (-1,1.7) node[anchor=north west] {(1,2)};
\draw (-1,3.15) node[anchor=north west] {(2,1)};
\draw [line width=2.pt] (5.,1.)-- (6.,0.);
\draw [line width=2.pt] (6.,0.)-- (8.,0.);
\draw [line width=2.pt] (8.,0.)-- (9.,1.);
\draw [line width=2.pt] (9.,1.)-- (9.,3.);
\draw [line width=2.pt] (9.,3.)-- (7.,4.);
\draw [line width=2.pt] (7.,4.)-- (5.,3.);
\draw [line width=2.pt] (5.,3.)-- (5.,1.);
\draw (6.3,4.85) node[anchor=north west] {(0,0)};
\draw (7.6,3.15) node[anchor=north west] {(1,1)};
\draw (7.6,1.7) node[anchor=north west] {(2,2)};
\draw (7.3,-0.1) node[anchor=north west] {(0,1)};
\draw (5,1.7) node[anchor=north west] {(1,2)};
\draw (5.3,-0.1) node[anchor=north west] {(2,0)};
\draw (5,3.15) node[anchor=north west] {(2,1)};
\draw [line width=2.pt] (11.,3.)-- (12.,4.);
\draw [line width=2.pt] (12.,4.)-- (14.,4.);
\draw [line width=2.pt] (14.,4.)-- (15.,3.);
\draw [line width=2.pt] (15.,3.)-- (15.,1.);
\draw [line width=2.pt] (15.,1.)-- (14.,0.);
\draw [line width=2.pt] (14.,0.)-- (12.,0.);
\draw [line width=2.pt] (12.,0.)-- (11.,1.);
\draw [line width=2.pt] (11.,1.)-- (11.,3.);
\draw (11.3,4.85) node[anchor=north west] {(0,0)};
\draw (13.3,4.85) node[anchor=north west] {(1,1)};
\draw (13.6,3.15) node[anchor=north west] {(2,2)};
\draw (13.6,1.7) node[anchor=north west] {(0,1)};
\draw (13.3,-0.1) node[anchor=north west] {(1,2)};
\draw (11.3,-0.1) node[anchor=north west] {(2,0)};
\draw (11,1.7) node[anchor=north west] {(0,2)};
\draw (11,3.15) node[anchor=north west] {(2,1)};
\begin{scriptsize}
\draw [fill=black] (1,0.) circle (2.5pt);
\draw [fill=black] (-1.,1.) circle (2.5pt);
\draw [fill=black] (-1.,3.) circle (2.5pt);
\draw [fill=black] (1,3.94) circle (2.5pt);
\draw [fill=black] (3.,1.) circle (2.5pt);
\draw [fill=black] (3.,3.) circle (2.5pt);
\draw [fill=black] (6.,0.) circle (2.5pt);
\draw [fill=black] (8.,0.) circle (2.5pt);
\draw [fill=black] (5.,1.) circle (2.5pt);
\draw [fill=black] (5.,3.) circle (2.5pt);
\draw [fill=black] (9.,1.) circle (2.5pt);
\draw [fill=black] (9.,3.) circle (2.5pt);
\draw [fill=black] (7.,4.) circle (2.5pt);
\draw [fill=black] (12.,0.) circle (2.5pt);
\draw [fill=black] (14.,0.) circle (2.5pt);
\draw [fill=black] (11.,1.) circle (2.5pt);
\draw [fill=black] (11.,3.) circle (2.5pt);
\draw [fill=black] (12.,4.) circle (2.5pt);
\draw [fill=black] (14.,4.) circle (2.5pt);
\draw [fill=black] (15.,3.) circle (2.5pt);
\draw [fill=black] (15.,1.) circle (2.5pt);
\end{scriptsize}
\end{tikzpicture}
\caption{Orthogonal Colourings of $C_6,C_7$, and $C_8$.}
\label{Figure: Case 1}
\end{figure}
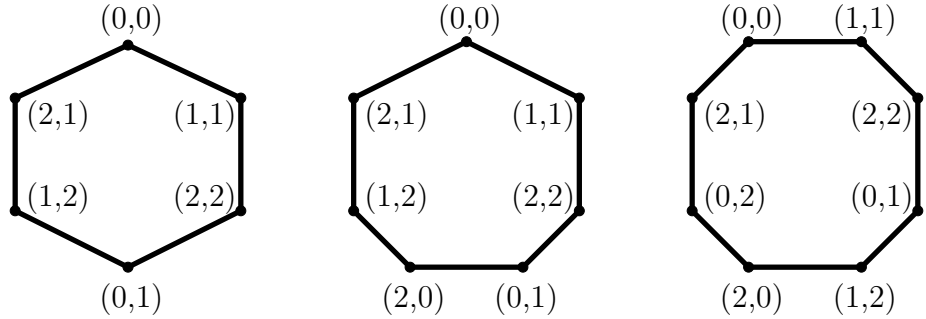

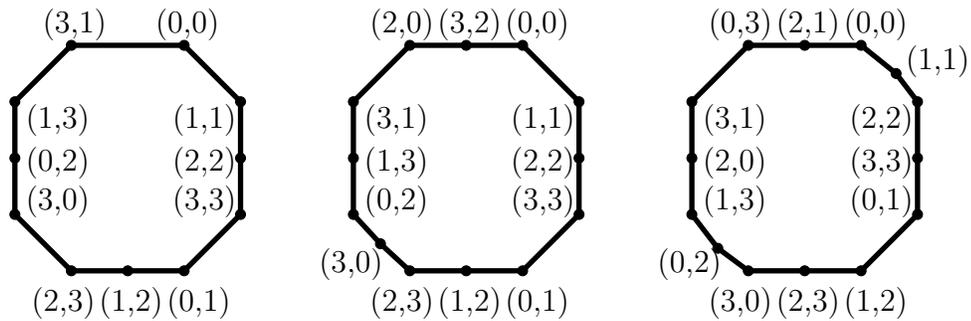
\begin{figure}[h!]
\centering
\begin{tikzpicture}[scale=.75,line cap=round,line join=round,>=triangle 45,x=1.0cm,y=1.0cm]
\draw [line width=2.pt] (0.,3.)-- (0.,2.);
\draw [line width=2.pt] (0.,2.)-- (0.,1.);
\draw [line width=2.pt] (0.,1.)-- (1.,0.);
\draw [line width=2.pt] (1.,0.)-- (2.,0.);
\draw [line width=2.pt] (2.,0.)-- (3.,0.);
\draw [line width=2.pt] (3.,0.)-- (4.,1.);
\draw [line width=2.pt] (4.,1.)-- (4.,2.);
\draw [line width=2.pt] (4.,2.)-- (4.,3.);
\draw [line width=2.pt] (0.,3.)-- (1.,4.);
\draw [line width=2.pt] (1.,4.)-- (3.,4.);
\draw [line width=2.pt] (3.,4.)-- (4.,3.);
\draw [line width=2.pt] (6.48,0.48)-- (6.,1.);
\draw [line width=2.pt] (6.,1.)-- (6.,2.);
\draw [line width=2.pt] (6.,2.)-- (6.,3.);
\draw [line width=2.pt] (6.,3.)-- (7.,4.);
\draw [line width=2.pt] (7.,4.)-- (8.,4.);
\draw [line width=2.pt] (8.,4.)-- (9.,4.);
\draw [line width=2.pt] (9.,4.)-- (10.,3.);
\draw [line width=2.pt] (10.,3.)-- (10.,2.);
\draw [line width=2.pt] (10.,2.)-- (10.,1.);
\draw [line width=2.pt] (10.,1.)-- (9.,0.);
\draw [line width=2.pt] (9.,0.)-- (8.,0.);
\draw [line width=2.pt] (8.,0.)-- (7.,0.);
\draw [line width=2.pt] (7.,0.)-- (6.48,0.48);
\draw [line width=2.pt] (12.46,0.4)-- (12.,1.);
\draw [line width=2.pt] (12.,1.)-- (12.,2.);
\draw [line width=2.pt] (12.,2.)-- (12.,3.);
\draw [line width=2.pt] (12.,3.)-- (13.,4.);
\draw [line width=2.pt] (13.,4.)-- (14.,4.);
\draw [line width=2.pt] (14.,4.)-- (15.,4.);
\draw [line width=2.pt] (15.,4.)-- (15.62,3.5);
\draw [line width=2.pt] (15.62,3.5)-- (16.,3.);
\draw [line width=2.pt] (16.,2.)-- (16.,3.);
\draw [line width=2.pt] (16.,2.)-- (16.,1.);
\draw [line width=2.pt] (16.,1.)-- (15.,0.);
\draw [line width=2.pt] (15.,0.)-- (14.,0.);
\draw [line width=2.pt] (14.,0.)-- (13.,0.);
\draw [line width=2.pt] (13.,0.)-- (12.46,0.4);
\draw (2.3,4.85) node[anchor=north west] {(0,0)};
\draw (2.6,3.15) node[anchor=north west] {(1,1)};
\draw (2.6,2.4) node[anchor=north west] {(2,2)};
\draw (2.6,1.7) node[anchor=north west] {(3,3)};
\draw (2.5,-0.1) node[anchor=north west] {(0,1)};
\draw (1.3,-0.1) node[anchor=north west] {(1,2)};
\draw (0.1,-0.1) node[anchor=north west] {(2,3)};
\draw (0,1.7) node[anchor=north west] {(3,0)};
\draw (0,2.4) node[anchor=north west] {(0,2)};
\draw (0,3.15) node[anchor=north west] {(1,3)};
\draw (.3,4.85) node[anchor=north west] {(3,1)};
\draw (8.5,4.85) node[anchor=north west] {(0,0)};
\draw (8.6,3.15) node[anchor=north west] {(1,1)};
\draw (8.6,2.4) node[anchor=north west] {(2,2)};
\draw (8.6,1.7) node[anchor=north west] {(3,3)};
\draw (8.5,-0.1) node[anchor=north west] {(0,1)};
\draw (7.3,-0.1) node[anchor=north west] {(1,2)};
\draw (6.1,-0.1) node[anchor=north west] {(2,3)};
\draw (5.2,0.6) node[anchor=north west] {(3,0)};
\draw (6,1.7) node[anchor=north west] {(0,2)};
\draw (6,2.4) node[anchor=north west] {(1,3)};
\draw (6.1,4.85) node[anchor=north west] {(2,0)};
\draw (6,3.15) node[anchor=north west] {(3,1)};
\draw (7.3,4.85) node[anchor=north west] {(3,2)};
\draw (14.5,4.85) node[anchor=north west] {(0,0)};
\draw (15.6,4.2) node[anchor=north west] {(1,1)};
\draw (14.6,3.15) node[anchor=north west] {(2,2)};
\draw (14.6,2.4) node[anchor=north west] {(3,3)};
\draw (14.6,1.7) node[anchor=north west] {(0,1)};
\draw (14.5,-0.1) node[anchor=north west] {(1,2)};
\draw (13.3,-0.1) node[anchor=north west] {(2,3)};
\draw (12,3.15) node[anchor=north west] {(3,1)};
\draw (12.1,-0.1) node[anchor=north west] {(3,0)};
\draw (11.2,0.6) node[anchor=north west] {(0,2)};
\draw (12,1.7) node[anchor=north west] {(1,3)};
\draw (12,2.4) node[anchor=north west] {(2,0)};
\draw (12.1,4.85) node[anchor=north west] {(0,3)};
\draw (13.3,4.85) node[anchor=north west] {(2,1)};
\begin{scriptsize}
\draw [fill=black] (1.,0.) circle (2.5pt);
\draw [fill=black] (0.,1.) circle (2.5pt);
\draw [fill=black] (0.,3.) circle (2.5pt);
\draw [fill=black] (1.,4.) circle (2.5pt);
\draw [fill=black] (3.,0.) circle (2.5pt);
\draw [fill=black] (4.,1.) circle (2.5pt);
\draw [fill=black] (4.,3.) circle (2.5pt);
\draw [fill=black] (0.,2.) circle (2.5pt);
\draw [fill=black] (4.,2.) circle (2.5pt);
\draw [fill=black] (2.,0.) circle (2.5pt);
\draw [fill=black] (3.,4.) circle (2.5pt);
\draw [fill=black] (8.,0.) circle (2.5pt);
\draw [fill=black] (7.,0.) circle (2.5pt);
\draw [fill=black] (6.,1.) circle (2.5pt);
\draw [fill=black] (6.,2.) circle (2.5pt);
\draw [fill=black] (6.,3.) circle (2.5pt);
\draw [fill=black] (7.,4.) circle (2.5pt);
\draw [fill=black] (8.,4.) circle (2.5pt);
\draw [fill=black] (9.,4.) circle (2.5pt);
\draw [fill=black] (10.,3.) circle (2.5pt);
\draw [fill=black] (10.,2.) circle (2.5pt);
\draw [fill=black] (10.,1.) circle (2.5pt);
\draw [fill=black] (9.,0.) circle (2.5pt);
\draw [fill=black] (6.48,0.48) circle (2.5pt);
\draw [fill=black] (12.,1.) circle (2.5pt);
\draw [fill=black] (12.,2.) circle (2.5pt);
\draw [fill=black] (12.,3.) circle (2.5pt);
\draw [fill=black] (13.,4.) circle (2.5pt);
\draw [fill=black] (14.,4.) circle (2.5pt);
\draw [fill=black] (15.,4.) circle (2.5pt);
\draw [fill=black] (16.,3.) circle (2.5pt);
\draw [fill=black] (16.,2.) circle (2.5pt);
\draw [fill=black] (16.,1.) circle (2.5pt);
\draw [fill=black] (15.,0.) circle (2.5pt);
\draw [fill=black] (14.,0.) circle (2.5pt);
\draw [fill=black] (13.,0.) circle (2.5pt);
\draw [fill=black] (12.46,0.4) circle (2.5pt);
\draw [fill=black] (15.62,3.5) circle (2.5pt);
\end{scriptsize}
\end{tikzpicture}
\caption{Orthogonal Colouring of $C_{11},C_{13}$, and $C_{14}$.}
\label{Figure: Case 2}
\end{figure}

For $n>4$, the remaining cases left to consider are $n=6,7,8,11,13,14$. These can be orthogonally coloured with $\lceil\sqrt{n}\,\rceil$ colours, shown in Figure \ref{Figure: Case 1} and Figure \ref{Figure: Case 2}. For $n=3,4$, $C_n$ needs $3$ colours. The following theorem summarizes these results.

\begin{thrm}
$O\chi(C_n)=\lceil\sqrt{n}\,\rceil$ if and only if $n>4$.
\end{thrm}
\begin{pf}
Follows from Lemma \ref{Lemma:Cycle Part 1}, Lemma \ref{Lemma: Cycle Part 2}, Figure \ref{Figure: Case 1}, and Figure \ref{Figure: Case 2}.
\end{pf}

Now that the orthogonal chromatic number of $C_n$ is determined, the goal is to construct multiple mutual orthogonal colourings. This is accomplished if restrictions are put upon the size of the vertex set. The following theorem gives a method to construct multiple mutual orthogonal colourings of $C_n$ when $\lceil\sqrt{n}\,\rceil=p$, where $p$ is a prime number.

\begin{thrm}
If $\lceil\sqrt{n}\,\rceil=p$ is a prime number, then $O\chi_{p-2}(C_n)=p$.
\end{thrm}
\begin{pf}
Let $v_i\in \mathbb{Z}_n$ where $0\leq i <n$. It is now shown that there are $p-1$ orthogonal assignments, but only $p-2$ of these assignments are proper colourings for a given $n$. Consider the following $p-1$ colourings: $$c_k(v_i)=\Big(i+k\Big\lfloor\frac{i}{p}\Big\rfloor\Big)(\textrm{mod}~p)$$ where $0\leq k < p$. It is now shown that any two are mutually orthogonal. Suppose otherwise, that is, $c_t(v_i)=c_t(v_j)$ and $c_s(v_i)=c_s(v_j)$, where $0\leq t<s<p$ and $0\leq i<j<n$. Then: 
\begin{align}
(i-j)&\equiv t\Big(\Big\lfloor\frac{j}{p}\Big\rfloor-\Big\lfloor\frac{i}{p}\Big\rfloor\Big)(\textrm{mod}~p). \label{Equation: 1}\\
(i-j)&\equiv s\Big(\Big\lfloor\frac{j}{p}\Big\rfloor-\Big\lfloor\frac{i}{p}\Big\rfloor\Big)(\textrm{mod}~p). \label{Equation: 2}
\end{align}

Note that $t,s\in \mathbb{Z}_p$, which is a field because $p$ is a prime. Thus, if $t\neq 0$, then $t^{-1}$ and $s^{-1}$ exist. Then, Equation \ref{Equation: 1} and Equation \ref{Equation: 2} can be rearranged to give that $(i-j)(s^{-1}-t^{-1})\equiv 0(\textrm{mod}~p)$. However, $0\leq t<s<p$, so $t^{-1}\not\equiv s^{-1}(\textrm{mod}~p)$. Therefore, it must be that $i\equiv j(\textrm{mod}~p)$, or equivalently, $j=i+mp$ where $0<m<p$. Since $s,m\neq 0$, $sm\not\equiv 0(\textrm{mod}~p)$ because $\mathbb{Z}_p$ has no zero divisors. This implies that:
\begin{align*}
\Big(j+s\Big\lfloor\frac{j}{p}\Big\rfloor\Big)(\textrm{mod}~p)&\equiv \Big(i+s\Big\lfloor\frac{i+mp}{p}\Big\rfloor\Big)(\textrm{mod}~p)\\
&\equiv \Big(i+sm+s\Big\lfloor\frac{i}{p}\Big\rfloor\Big)(\textrm{mod}~p)\\
&\not\equiv \Big(i+s\Big\lfloor\frac{i}{p}\Big\rfloor\Big)(\textrm{mod}~p).
\end{align*}

This contradicts Equation \ref{Equation: 2} however, hence $t=0$. Putting $t=0$ into Equation \ref{Equation: 1} gives that $i\equiv j(\textrm{mod}~p)$, and the same contradiction arises. Therefore, the colourings are all mutually orthogonal. It remains to show that $p-2$ of the colourings are proper.

Notice that $k\Big\lfloor\frac{i}{p}\Big\rfloor\leq k\Big\lfloor\frac{i+1}{p}\Big\rfloor\leq k\Big\lfloor\frac{i}{p}\Big\rfloor+k $. Therefore, $1\leq c_k(v_{i+1})-c_k(v_i)\leq k$ where $1\leq i \leq n-2$. Since $k<p$, this gives that $c_k(v_{i+1})\neq c_k(v_i)$. It remains to show that there is no colour conflict between the vertices  $v_0$ and $v_{n-1}$. Note that the colour $0$ is assigned to $v_0$ in all of the colourings. Therefore, by the mutual orthogonality of the colourings, at most one colouring has $0$ assigned to the vertex $v_{n-1}$. Therefore, choosing the $k-1=p-2$ colourings that don't assign the colour $0$ to the vertex $v_{n-1}$ gives a $(p-2)$-orthogonal colouring of $C_n$ using $p$ colours.
\end{pf}

\section{Orthogonal Colourings of Circulant Graphs}

Cayley graphs on the additive group $\mathbb{Z}_{p^2}$ are now considered, where $p$ is a prime. Also, the size of the generating set $S$ is varied throughout this section. To start, a method for orthogonally colouring $\Gamma(\mathbb{Z}_{p^2},S)$ is constructed when $|S|$ is sufficiently small. Let $\alpha\in\mathbb{Z}_p\backslash\{0\}$ and consider the function $\hat F_{\alpha,p}:\mathbb{Z}_p\times \mathbb{Z}_p\to \mathbb{Z}_{p^2}$ defined by $$\hat F_{\alpha,p}(i,j)=(((\alpha(j-i)(\textrm{mod}~p))+p(2i-j))(\textrm{mod}~p^2).$$ This function assigns colour pairs from $\mathbb{Z}_p\times\mathbb{Z}_p$ to the vertices of $\Gamma(\mathbb{Z}_{p^2},S)$. Therefore, the goal is to show that the inverse of this function is an orthogonal colouring. The following lemma shows that this assignment is injective.

\begin{lem}\label{Lemma: Circulant Map 2}
For every $\alpha\in\mathbb{Z}_p\backslash\{0\}$, $\hat F_{\alpha,p}(i,j)$ is a bijection.
\end{lem}

\begin{pf}
Since $|\mathbb{Z}_p\times \mathbb{Z}_p|=|\mathbb{Z}_{p^2}|=p^2$, it is sufficient to show that $\hat F_{\alpha,p}$ is an injective function. Suppose  $\hat F_{\alpha,p}(i,j)=\hat F_{\alpha,p}(r,s)$. For this equality to be true, the two modular components of $\hat F_{\alpha,p}$ must be equal. That is:
\begin{align}
\alpha(j-i)(\textrm{mod}~p)&\equiv \alpha(s-r)(\textrm{mod}~p). \label{Equation: 3}\\
p(2i-j)(\textrm{mod}~p^2)&\equiv p(2r-s)(\textrm{mod}~p^2). \label{Equation: 4}
\end{align}

Note that $\alpha$ has a multiplicative inverse in $\mathbb{Z}_p$, which is denoted by $\alpha^{-1}$. Therefore, Equation (\ref{Equation: 3}) can be rewritten as $(i-j)(\textrm{mod}~p)\equiv (r-s)(\textrm{mod}~p)$. Multiplying by $p$ gives that $p(i-j)(\textrm{mod}~p^2)\equiv p(r-s)(\textrm{mod}~p^2)$. Substituting this into Equation (\ref{Equation: 4}) gives that $pi(\textrm{mod}~p^2)\equiv pr(\textrm{mod}~p^2)$. Note that this implies that $i(\textrm{mod}~p)\equiv r(\textrm{mod}~p)$. Since $i$ and $r$ are elements of $\mathbb{Z}_p$, it follows that $i=r$. Lastly, substituting $i=r$ into $(i-j)(\textrm{mod}~p)\equiv (r-s)(\textrm{mod}~p)$ gives that $j(\textrm{mod}~p)=s(\textrm{mod}~p)$. Similarly, since $j$ and $s$ are elements of $\mathbb{Z}_p$, it follows that $j=s$. Therefore, $(i,j)=(r,s)$.
\end{pf}

In particular, $\hat F_{\alpha,p}^{-1}:\mathbb{Z}_{p^2} \to \mathbb{Z}_p\times \mathbb{Z}_p$ is an injective function. Since $\hat F_{\alpha,p}^{-1}$ is injective into $\mathbb{Z}_p\times\mathbb{Z}_p$, $\hat F_{\alpha,p}^{-1}$ will orthogonally assign the vertices to colour pairs in $\mathbb{Z}_p\times \mathbb{Z}_p$. Therefore, if it can be shown that the components of $\hat F_{\alpha,p}^{-1}$ are both proper colourings, then $\hat F_{\alpha,p}^{-1}$ will be an orthogonal colouring. For a general generating set $S$ and $\alpha$, this is not always the case. However, the following theorem shows that if the size of the generating set $S$ is sufficiently small, then there is an $\alpha$ for which $\hat F_{\alpha,p}^{-1}$ is a proper colouring of $\Gamma(\mathbb{Z}_{p^2},S)$. 

\begin{thrm}\label{Theorem: Circulant 1}
For $p$ a prime, if $|S|<\frac{p-1}{2}$, then $O\chi(\Gamma(\mathbb{Z}_{p^2},S))=p$.
\end{thrm}
\begin{pf}
The goal is to show that $\hat F_{\alpha,p}^{-1}$ is an orthogonal colouring for some $\alpha\in \mathbb{Z}_{p}\backslash\{0\}$. By Lemma \ref{Lemma: Circulant Map 2}, $\hat F_{\alpha,p}^{-1}$ is an orthogonal assignment of the vertices. It remains to show that there is an $\alpha$ such that $\hat F_{\alpha,p}^{-1}$ is a proper colouring. Suppose two vertices $k$ and $l$ receive the same colour in the first colouring. That is, $\hat F_{\alpha,p}^{-1}(k)=(i,(j+x)(\textrm{mod}~p))$ and $\hat F_{\alpha,p}^{-1}(l)=(i,j)$ for some $i,j\in \mathbb{Z}_p$ and $x\in\mathbb{Z}_p\backslash\{0\}$. Then $$k-l=\hat F_{\alpha,p}(i,(j+x)(\textrm{mod}~p))-\hat F_{\alpha,p}(i,j)=(((\alpha x)(\textrm{mod}~p))-px)(\textrm{mod}~p^2).$$ 

Similarly, suppose two vertices $k$ and $l$ receive the same colour in the second colouring. That is, $\hat F_{\alpha,p}^{-1}(k)=((i+x)(\textrm{mod}~p),j)$ and $\hat F_{\alpha,p}^{-1}(l)=(i,j)$ for some $i,j\in \mathbb{Z}_p$ and $x\in\mathbb{Z}_p\backslash\{0\}$. Then $$k-l=\hat F_{\alpha,p}((i+x)(\textrm{mod}~p),j)-\hat F_{\alpha,p}(i,j)=(((-\alpha x)(\textrm{mod}~p))+2px)(\textrm{mod}~p^2).$$ Let $A_\alpha$ be the set of differences in the first colouring and $B_\alpha$ be the set of differences in the second colouring. That is, $A_{\alpha}=\{(((\alpha x)(\textrm{mod}~p))-px)(\textrm{mod}~p^2)|x\in\mathbb{Z}_p\backslash\{0\}\}$ and $B_{\alpha}=\{(((-\alpha x)(\textrm{mod}~p))+2px)(\textrm{mod}~p^2)|x\in\mathbb{Z}_p\backslash\{0\}\}$. Therefore, there is a colour conflict in the first colouring if and only if $S\cap A_\alpha\neq \emptyset$. Similarly, there is a colour conflict in the second colouring in and only if $S\cap B_\alpha \neq \emptyset$.
 
Properties of $A_{\alpha}$ and $B_{\alpha}$ are now discussed. The first is that $A_\alpha = B_{2\alpha(\textrm{mod}~p)}$. Since $\mathbb{Z}_p$ is a field, $2$ has a multiplicative inverse in $\mathbb{Z}_p$, which is denoted $2^{-1}$. Then, since $\mathbb{Z}_p$ has no zero divisors, $\mathbb{Z}_p\backslash\{0\}=-2\mathbb{Z}_p\backslash\{0\}$. Therefore, it follows that
\begin{align*}
B_{2\alpha(\textrm{mod}~p)}&=\{(((-2\alpha x)(\textrm{mod}~p))+2px)(\textrm{mod}~p^2)|x\in\mathbb{Z}_p\backslash\{0\}\}\\
&=\{(((-2\alpha (-2^{-1}x)(\textrm{mod}~p))+2p(-2^{-1}x)(\textrm{mod}~p^2)|x\in -2\mathbb{Z}_p\backslash\{0\}\}\\
&=\{(((\alpha x)(\textrm{mod}~p))-px)(\textrm{mod}~p^2)|x\in -2\mathbb{Z}_p\backslash\{0\}\}\\
&=\{(((\alpha x)(\textrm{mod}~p))-px)(\textrm{mod}~p^2)|x\in\mathbb{Z}_p\backslash\{0\}\}\\
&=A_\alpha.
\end{align*}

A second property of $A_\alpha$ is the following. The $A_\alpha$'s together with $\{m| m\in\mathbb{Z}_p\backslash\{0\}\}$ and $\{mp| m\in\mathbb{Z}_p\backslash\{0\}\}$ are disjoint. First, notice that $\alpha x(\textrm{mod}~p)\not\equiv 0(\textrm{mod}~p)$ for any $x\in\mathbb{Z}_p\backslash\{0\}$. Therefore, it follows that $A_{\alpha}\cap \{mp|m\in\mathbb{Z}_p\backslash\{0\}\}=\emptyset$. Next, notice that $-px(\textrm{mod}~p^2)\not\equiv 0(\textrm{mod}~p^2)$ for any $x\in\mathbb{Z}_p\backslash\{0\}$. Therefore, $A_\alpha\cap\{m| m\in\mathbb{Z}_p\backslash\{0\}\}=\emptyset$. It is also the case that $\{mp| m\in\mathbb{Z}_p\backslash\{0\}\}\cap\{m| m\in\mathbb{Z}_p\backslash\{0\}\}=\emptyset$. This is because the first set is multiples of $p$ and the second set is not. Therefore, it remains to show that the $A_\alpha$'s are all mutually disjoint.

Let $\alpha_1,\alpha_2\in \mathbb{Z}_p\backslash\{0\}$ where $\alpha_1\neq \alpha_2$. Suppose that $A_{\alpha_1}\cap A_{\alpha_2}\neq\emptyset$. Then, there exists some $c=((\alpha_1 x(\textrm{mod}~p))-px)(\textrm{mod}~p^2)$ and $c=((\alpha_2 y(\textrm{mod}~p))-py)(\textrm{mod}~p^2)$ for some $x,y\in\mathbb{Z}_p\backslash\{0\}$. Note that $c$ can be uniquely written as $r-px$, where $r\in\mathbb{Z}_p\backslash\{0\}$. This means that $-px(\textrm{mod}~p^2)=-py(\textrm{mod}~p^2)$. This implies that $x(\textrm{mod}~p)=y(\textrm{mod}~p)$. Since $x$ and $y$ are elements of $\mathbb{Z}_p$, its follows that $x=y$. However, since $c$ can be uniquely written as $r-px$, it follows that $\alpha_1 x(\textrm{mod}~p)=\alpha_2 y(\textrm{mod}~p)$. Substituting $x=y$ into this gives that $\alpha_1(\textrm{mod}~p)=\alpha_2(\textrm{mod}~p)$. Then, since $\alpha_1$ and $\alpha_2$ are elements of $\mathbb{Z}_p\backslash\{0\}$, it follows that $\alpha_1=\alpha_2$, a contradiction. Therefore, the $A_\alpha$'s along with $\{mp| m\in\mathbb{Z}_p\backslash\{0\}\}$ and $\{m| m\in\mathbb{Z}_p\backslash\{0\}\}$ are all mutually disjoint.

Now to show that there is a choice of $\alpha$ for which $\hat F_{\alpha,p}^{-1}$ is a proper colouring. Recall that there is a colour conflict in the first colouring if and only if $S\cap A_\alpha\neq \emptyset$ and there is a colour conflict in the second colouring if and only if $S\cap B_{\alpha}\neq \emptyset$. Let $c\in S$ and suppose that $c\in \{mp| m\in\mathbb{Z}_p\backslash\{0\}\}$ or $c\in\{m| m\in\mathbb{Z}_p\backslash\{0\}\}$. Then, since the $A_\alpha$'s, along with $\{mp| m\in\mathbb{Z}_p\backslash\{0\}\}$ and $\{m| m\in\mathbb{Z}_p\backslash\{0\}\}$ are all mutually disjoint, $c\not\in A_\alpha$ for any $\alpha\in\mathbb{Z}_p\backslash\{0\}$. Since $B_{2\alpha(\textrm{mod}~p)}=A_\alpha$, then $c\not\in B_\alpha$ for any $\alpha$ as well. Therefore, $c$ will not result in any colouring conflicts in the first colouring or in the second colouring.

Suppose now that $c\in A_{\alpha_1}$ for some $\alpha_1\in\mathbb{Z}_p\backslash\{0\}$. Then, since the $A_\alpha$'s are disjoint, for any $\alpha_2\neq \alpha_1$, $c\not\in A_{\alpha_2}$. Thus, since $B_{2\alpha(\textrm{mod}~p)}=A_{\alpha}$, it follows that $c\in B_{2\alpha_1}$ but $c\not\in B_{\alpha_3}$ for any $\alpha_3\neq 2\alpha_1$. Hence, if $\alpha_1$ is chosen, then $c\in A_{\alpha_1}$, and there is a colour conflict in the first colouring. Similarly, if $2\alpha_1$ is chosen, then $c\in B_{2\alpha_1}$, and there is a colour conflict in the second colouring. However, any other choice of $\alpha$ will result in $c\not\in A_\alpha$ and $c\not\in B_\alpha$. Thus each $c\in S$ will result in at most $2$ restrictions on the choice of $\alpha$. Since that are at most $|S|<\frac{p-1}{2}$ elements in $S$, there are fewer than $p-1$ restrictions on the choice of $\alpha$. Since $\alpha\in\mathbb{Z}_p\backslash\{0\}$, there are $p-1$ choices for $\alpha$. Therefore, there are more choices than restrictions. 
\end{pf}

Theorem \ref{Theorem: Circulant 1} says that if the size of the generating set is sufficiently small, then an orthogonal colouring can be constructed using only $p$ colours. This leads to the following question.

\begin{ques}
What is the largest $m$, such that if $|S|< m$, then $O\chi(\Gamma(\mathbb{Z}_{p^2},S))=p$?
\end{ques}
  Theorem \ref{Theorem: Circulant 1} gives a lower bound of $m> \frac{p-1}{2}$. However, on the other hand, the following theorem by Klotz and Sander gives an upper bound on the value of $m$.

\begin{thrm}[Klotz and Sander \cite{klotz2016uniquely}]\label{Theorem: Walter}
If $S=\{\pm 1,\pm 2,\dots, \pm (p-1)\}$, then $\Gamma(\mathbb{Z}_{p^2},S)$ is uniquely $p$-colourable.
\end{thrm}

Note that uniquely $p$-colourable graphs can not be orthogonally coloured with $p$ colours, unless the graph is $K_p$. This is because the first and second colouring are the same, and thus if any colour is used twice in the first colouring, then that colour pair occurs twice. Therefore, since $|S|$ in Theorem \ref{Theorem: Walter} is of size $2p-2$, this gives an upper bound of $m<2p-2$.

It is now shown that if no multiples of $p$ are in $S$, then the lower bound on $m$ can be increased to $p-1$. Let $\alpha\in\mathbb{Z}_{p^2}$ be such that $gcd(p^2,\alpha)=1$. Then, consider the following function: $F_{\alpha,p}:\mathbb{Z}_p\times \mathbb{Z}_p\to \mathbb{Z}_{p^2}$ defined by $$F_{\alpha,p}(i,j)=(ip+j\alpha)(\textrm{mod}~p^2).$$ This function will be used to assign colour pairs to $\Gamma(\mathbb{Z}_{p^2},S)$, when there are no multiples of $p$ in $S$. The following lemma shows that this assignment is injective.

\begin{lem}\label{Lemma: Circulant Map 1}
For any $\alpha\in\mathbb{Z}_{p^2}\backslash\{xp|0\leq x<p\}$, $F_{\alpha,p}(i,j)$ is a bijection.
\end{lem}
\begin{pf}
Since $|\mathbb{Z}_p\times \mathbb{Z}_p|=|\mathbb{Z}_{p^2}|=p^2$, it is sufficient to show that $F_{\alpha,p}$ is surjective. Let $x\in \mathbb{Z}_{p^2}$. Since $gcd(p^2,\alpha)=1$, the multiplicative inverse $\alpha^{-1}\in \mathbb{Z}_{p^2}$. Now, by the Division Algorithm, there exist unique integers $q$ and $r$ such that $x=qp+r$ where $0\leq q,r <p$.  Let $i=q$ and $j=r\alpha^{-1}(\textrm{mod}~p^2)$. Substituting gives that $F_{\alpha,p}(i,j)=(qp+r)(\textrm{mod}~p^2)=x(\textrm{mod}~p^2)$. Therefore, $F_{\alpha,p}$ is surjective, and thus bijective.
\end{pf}

In particular, $F_{\alpha,p}^{-1}:\mathbb{Z}_{p^2} \to \mathbb{Z}_n\times \mathbb{Z}_n$ is an injective function. Therefore, $F_{\alpha,p}^{-1}$  will orthogonally assign the vertices to colour pairs in $\mathbb{Z}_p\times\mathbb{Z}_p$. The following lemma shows that if the size of the generating set $S$ is sufficiently small and $S$ contains no multiples of $p$, then there is an $\alpha$ for which $F_{\alpha,p}^{-1}$ is a proper colouring of $\Gamma(\mathbb{Z}_{p^2},S)$.

\begin{thrm}\label{Theorem: Circulant 1}
For $p$ a prime, if $|S|<p$ and $xp\not\in S$ for all $x$ where $x\in \mathbb{Z}_p\backslash\{0\}$, then $O\chi(\Gamma(\mathbb{Z}_{p^2},S))=p$.
\end{thrm}
\begin{pf}
It is shown that $F_{\alpha,p}^{-1}$ is an orthogonal colouring for some $\alpha\in \mathbb{Z}_{p^2}$ where $\gcd(p^2,\alpha)=1$. By Lemma \ref{Lemma: Circulant Map 1}, $F_{\alpha,p}^{-1}$ is an orthogonal assignment of the vertices. It remains to show that there is a choice for $\alpha$ such that $F_{\alpha,p}^{-1}$ is a proper colouring. 

Suppose two vertices $k$ and $l$ receive the same colour in the first colouring. That is, $F^{-1}_{\alpha,p}(k)=(i,(j+x)(\textrm{mod}~p))$ and $F^{-1}_{\alpha,p}(l)=(i,j)$ for some $i,j\in \mathbb{Z}_p$ and $x\in \mathbb{Z}_p\backslash\{0\}$. Note that $F_{\alpha,p}(i,(j+x)(\textrm{mod}~p))=ip+(j+x)\alpha=k$ and $F_{\alpha,p}(i,j)=ip+j\alpha=l$. Therefore, there is a colour conflict in the first colouring if and only if $k$ and $l$ are adjacent, which occurs when $k-l=x\alpha \in S$.

Suppose two vertices $k$ and $l$ receive the same colour in the second colouring. That is, $F^{-1}_{\alpha,p}(k)=(i,j)$ and $F^{-1}_{\alpha,p}(l)=((i+x)(\textrm{mod}~p),j)$ for some $i,j\in \mathbb{Z}_p$ and $x\in\mathbb{Z}_p\backslash\{0\}$. Note that $F_{\alpha,p}(i,j)=ip+j\alpha=k$ and $F_{\alpha,p}((i+x)(\textrm{mod}~p),j)=(i+x)p+j\alpha=l$. Therefore, there is a colour conflict in the second colouring if and only if $k$ and $l$ are adjacent, so when $l-k=xp \in S$. By assumption, $xp\not\in S$ for all $x$, so there are no colour conflicts in the second colouring.

So the only conflict that can occur is in the first colouring, which happens when $x\alpha \in S$. Since $x\in \mathbb{Z}_p\backslash\{0\}$, which is a field, $x^{-1}$ exists. Therefore, if $\alpha\not\in\bigcup_{x}\{x^{-1}S\}$, then there will be no colour conflicts. Note that since $p$ is a prime, there are $p(p-1)$ choices for $\alpha$ so that $\gcd(p^2,\alpha)=1$. Then, since there are $p-1$ choices for $x$, there are at most $(p-1)|S|<(p-1)p$ elements in $\bigcup_x\{x^{-1}S\}$.  Therefore, there is an $\alpha$ such that $\alpha\not\in\bigcup_{x}\{x^{-1}S\}$. 
\end{pf}

This concludes our study of circulant graphs with varying generating sets. Certain circulant graphs that have large generating sets are now studied. This is illustrated with the Paley graphs. To talk about this graph, relevant properties of finite fields are now summarized.
\section{Paley Graphs}

Finite fields of order $q$ exist if and only if $q=p^k$ where $p$ is a prime and $k\in\mathbb{Z}^+$. These fields are unique up to isomorphism, so they are denoted  $\mathbb{F}_q$. The multiplicative group, $\mathbb{F}_q^*$, is cyclic, so all non-zero elements can be expressed as powers of a single element, called a primitive element of the field.

Finite fields can be explicitly constructed as such. If $q=p^k$, then $\mathbb{F}_q\cong \mathbb{Z}_p[X]/(P)$ where $(P)$ is the ideal generated by an irreducible polynomial $P$ of degree $k$ in $\mathbb{Z}_p[X]$. The Paley graph, denoted $QR(q)$, can be constructed as a Cayley graph of $\mathbb{F}_q$. Let $\alpha$ be a primitive element of $\mathbb{F}_q$ and let $S=\{\alpha^{2m}:1\leq m\leq \frac{q-1}{2}\}$. That is, $S$ is the set of all quadratic residues in $\mathbb{F}_q$. Then the Paley graph $QR(q)=\Gamma(\mathbb{F}_q,S)$. 
 
Paley graphs have a variety of interesting properties, but the one that is used is that they are self-complementary \cite{jones2016paley}. This allows the cliques in the Paley graphs to be considered as colour classes, since they can be turned into independent sets by taking the complement. The following lemma, which is an exercise in \cite{fraleigh2003first}, describes a relation between cosets of two different subgroups.
 
\begin{lem}\label{Lemma: Coset Intersection}
Let $G$ be a group, $H\leq G$, and $K\leq G$. Then for any $a,b,c\in G$, either $(a+H)\cap (b+K)=\emptyset$ or $(a+H)\cap (b+K)= c+(H\cap K)$.
\end{lem}

The following theorem gives a method for constructing multiple orthogonal colourings of Paley graphs by utilizing the structure of the finite field. The idea is to find particular cosets that share at most one element. Then, Lemma \ref{Lemma: Coset Intersection} gives that their intersection also contains at most one element. These cosets can then be used as the colour classes. Since they intersect in at most one element, they will be an orthogonal assignment of the colours. The following theorem formalizes this argument.

\begin{thrm}
For $p$ a prime and an integer $r\geq 1$, $O\chi_{\frac{p^r+1}{2}}(QR(p^{2r}))=p^r$.
\end{thrm}
\begin{pf}
Let $G=\mathbb{F}_{p^{2r}}$. Since $p^r\mid p^{2r}$, there exists a subfield $H\subset G$ where $H\cong \mathbb{F}_{p^r}$. Also, notice that $H^*$ is a multiplicative subgroup of $G^*$ with index $\frac{p^{2r}-1}{p^r-1}=p^r+1$. Thus, $H=\{0\}\bigcup\{\alpha^{m(p^r+1)}|1\leq m < p^r-1\}$ where $\alpha$ is a primitive element of $G$. Next, for $0\leq i < \frac{p^r+1}{2}$, consider the sets $H_i=\{0\}\bigcup\{\alpha^{m(p^r+1)+2i}|1\leq m < p^r-1\}$. The sets $H_i$ are additive subgroups of $G$. To see this, it suffices to show for $x,y\in H_i$ that $x-y\in H_i$. \begin{align*}
x-y&=\alpha^{m(p^r+1)+2i}-\alpha^{n(p^r+1)+2i}  &&\text{by the definition of $H_i$}\\
&=\alpha^{2i}(\alpha^{m(p^r+1)}-\alpha^{n(p^r+1)})\\
&=\alpha^{2i}(\alpha^{t(p^r+1)})  &&\text{because $\mathbb{F}_{p^{2r}}$ is a field}\\
&=\alpha^{t(p^r+1)+2i}\\
&\in H_i &&\text{by the definition of $H_i$.}
\end{align*}

Therefore, the sets $H_i$ are additive subgroups of $G$. Since $\alpha^{m(p^r+1)+2i}$ is an even power of $\alpha$, $H_i\subset S$ and so $H_i$ are cliques in $QR(p^{2r})$. Similarly, the $p^r-1$ cosets of $H_i$ are cliques in $QR(p^{2r})$. Therefore, $H_i$ and its cosets are independent sets in the complement graph. For $1\leq j \leq p^r-1$, let $H_i+k_j$ be the $j$th coset of $H_i$. Then, define the colourings of the complement as $c_i(x)=j$ for $x\in H_i+k_j$. Since these are independent sets, the $c_i$ are proper colourings. It remains to show that these colourings are mutually orthogonal.

Notice that $H_i\bigcap H_j=\{0\}$. Therefore, by Lemma \ref{Lemma: Coset Intersection}, any coset of $H_i$ and any coset of $H_j$ will intersect in at most 1 element. Therefore, the colour pair $(i,j)$ is only assigned to at most one vertex. Hence, there is no orthogonal conflict.
\end{pf}

\section{Orthogonal Colourings of Product Graphs}

The Cartesian product of two graphs $G$ and $H$, $G\square H$, has vertex set $V(G)\times V(H)$, and two vertices $(u_1,v_1)$ and $(u_2,v_2)$ in $G\square H$ are adjacent if and only if either $u_1=u_2$ and $v_1v_2\in E(H)$ or $v_1=v_2$ and $u_1u_2\in E(G)$. A relationship between Cayley graphs and the Cartesian graph product was given by Theron \cite{theron1989extension}. He showed that if $G_1$ and $G_2$ are groups with Cayley graphs $\Gamma(G_1,S_1)$ and $\Gamma(G_2,S_2)$ respectively, then a Cayley graph of $G_1\times G_2$ is $\Gamma(G_1\times G_2,(S_1\times\{0\})\cup(\{0\}\times S_2))=\Gamma(G_1,S_1)\square \Gamma(G_2,S_2)$.

Let $S=\mathbb{Z}_q\backslash\{0\}$ and consider the group $\mathbb{Z}_q^d$ with addition (component-wise) as the group operation. Notice that $\Gamma(\mathbb{Z}_q,S)=K_q$. Therefore, by the above result, a possible Cayley graph for $\mathbb{Z}_q^d$ is the Cartesian graph product of $d$ complete graphs $K_q$. This is what is defined to be the Hamming graph $H(d,q)$. The following theorem shows how  the Cartesian graph product affects orthogonal colourings.

\begin{thrm}\label{Theorem: Cartesian Product}
If $G$ has $n^2$ vertices, $H$ has $m^2$ vertices, and $O\chi(G)=n\geq m$, then $O\chi(G\square H)=nm$.
\end{thrm}
\begin{pf}
Label $V(G)=\{v_k:0\leq k <n^2\}$ and $V(H)=\{(u_i,u_j):0\leq i,j<m \}$ where $n\geq m$. Let $f=(f_1,f_2)$ be an orthogonal colouring of $G$ using the colours $0,1,\dots, n-1$. It is shown that $g=(g_1,g_2)$ is an orthogonal colouring of $G\square H$ where:
\begin{align*}
g_1((v_k,(u_i,u_j)))&=(f_1(v_k)+j)(\textrm{mod}~n)+in.\\
g_2((v_k,(u_i,u_j)))&=(f_2(v_k)+i)(\textrm{mod}~n)+jn.
\end{align*}
First, it is shown that $g$ has no orthogonal conflicts. Let $v_{k_1},v_{k_2}\in V(G)$ and let $(u_{i_1},u_{j_1}),(u_{i_2},u_{j_2})\in V(H)$. If $g((v_{k_1},(u_{i_1},u_{j_1})))=g((v_{k_2},(u_{i_2},u_{j_2})))$, then:
\begin{align}
(f_1(v_{k_1})+j_1)(\textrm{mod}~n)+i_1n&=(f_1(v_{k_2})+j_2)(\textrm{mod}~n)+i_2n. \label{Equation: ONE}\\
(f_2(v_{k_1})+i_1)(\textrm{mod}~n)+j_1n&=(f_2(v_{k_2})+i_2)(\textrm{mod}~n)+j_2n. \label{Equation: TWO}
\end{align}
Without loss of generality, suppose that $i_1<i_2$. Then:
\begin{align*}
(f_1(v_{k_1})+j_1)(\textrm{mod}~n)+i_1n &<n+i_1n\\
&=n(1+i_1)\\
&\leq i_2 n\\
&\leq (f_1(v_{k_2})+j_2)(\textrm{mod}~n)+i_2n.
\end{align*}
Therefore, by transitivity, $(f_1(v_{k_1})+j_1)(\textrm{mod}~n)+i_1n<(f_1(v_{k_2})+j_2)(\textrm{mod}~n)+i_2n$, which contradicts Equation \eqref{Equation: ONE}, thus $i_1=i_2$. A similar argument shows that $j_1=j_2$. Substituting $i_1=i_2$ and $j_1=j_2$ into Equations \eqref{Equation: ONE} and \eqref{Equation: TWO}, gives $f_1(v_{k_1})=f_1(v_{k_2})$ and $f_2(v_{k_1})=f_2(v_{k_2})$. Hence, $v_{k_1}=v_{k_2}$ because $f$ is an orthogonal colouring. Thus, $(v_{k_1},(u_{i_1},u_{j_1}))=(v_{k_2},(u_{i_2},u_{j_2}))$.

It remains to show that $g_1$ and $g_2$ are proper colourings of $G\square H$. First, consider adjacencies of the form $(v_{k_1},(u_i,u_j))\sim(v_{k_2},(u_i,u_j))$. Since $v_{k_1}\sim v_{k_2}$ in $G$ and $f_1$ is a proper colouring of $G$, $f_1(v_{k_1})\neq f_1(v_{k_2})$. Thus 
\begin{align*}
g_1((v_{k_1},(u_i,u_j)))&=(f_1(v_{k_1})+j)(\textrm{mod}~n)+in \\
&\neq (f_1(v_{k_2})+j)(\textrm{mod}~n)+in \\
&= g_1((v_{k_2},(u_i,u_j)).
\end{align*}

Next, consider adjacencies of the form $(v_k,(u_{i_1},u_{j_1}))\sim (v_k,(u_{i_2},u_{j_2}))$. Suppose that $i_1=i_2$ and $j_1\neq j_2$, then:
\begin{align*}
g_1((v_k,(u_{i_1},u_{j_1})))&=(f_1(v_k)+ j_1)(\textrm{mod}~n) +i_1n\\
&\neq (f_1(v_k)+ j_2)(\textrm{mod}~n) +i_1n   &&\text{because $m\leq n$}\\
&=(f_1(v_k)+ j_2)(\textrm{mod}~n) +i_2n\\
&=g_1((v_k,(u_{i_2},u_{j_2}))).
\end{align*}

Thus, there are no colour conflicts in the case. If $i_1\neq i_2$, then the argument used to prove the orthogonality of $g$ shows that $g_1((v_k,(u_{i_1},u_{j_1})))\neq g_1((v_k,(u_{i_2},u_{j_2})))$. Hence, there are no conflicts in this case either. Therefore, $O\chi(G\square H)=nm$.
\end{pf}

Note that $H(d+2,q)=H(d,q)\square H(2,q)$. Therefore, if $O\chi(H(2,q))=q$, then Theorem \ref{Theorem: Cartesian Product} allows the use of induction to orthogonally colour Hamming graphs of the form $H(2d,q)$. Notice that a pair of orthogonal Latin squares corresponds to an orthogonal colouring of $K_q\square K_q= H(2,q)$. It is known that orthogonal Latin squares of size $n\neq 2,6$ exist \cite{bose1960further,tarry1900probleme}, giving the following lemma.

\begin{lem}\label{Lemma: Hamming Graph}
If $q\neq 2,6$ then $O\chi(H(2,q))=q$.
\end{lem}

\begin{cor}\label{Corollary: Case 1}
If $q\neq 2,6$, and $d\geq 1$, then $O\chi(H(2d,q))=q^d$.
\end{cor}
\begin{pf}
Proceed by induction on $d$. If $d=1$, then $O\chi(H(2,q))=q$ by Lemma \ref{Lemma: Hamming Graph}. Assume true for $d\leq k$, $k\geq 1$. Then, $H(2(k+1),q)=H(2k,q)\square H(2,q)$. By the induction hypothesis, $O\chi(H(2k,q))=q^k$. Therefore, by Theorem \ref{Theorem: Cartesian Product} it follows that, $O\chi(H(2(k+1),q)=q^{k+1}$.
\end{pf}

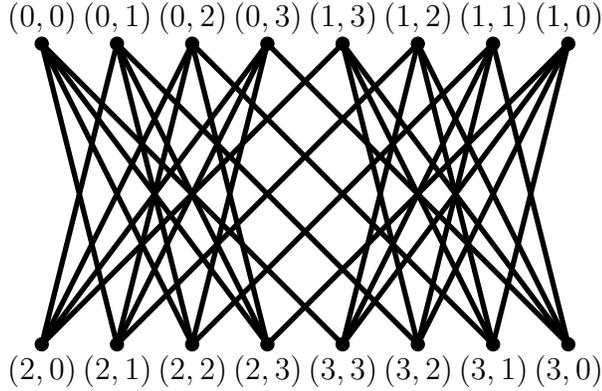
\begin{figure}[h!]
\centering
\begin{tikzpicture}[line cap=round,line join=round,>=triangle 45,x=1cm,y=1cm]
\draw [line width=2pt] (0,4)-- (1,0);
\draw [line width=2pt] (0,4)-- (2,0);
\draw [line width=2pt] (0,4)-- (3,0);
\draw [line width=2pt] (1,4)-- (0,0);
\draw [line width=2pt] (1,4)-- (2,0);
\draw [line width=2pt] (1,4)-- (3,0);
\draw [line width=2pt] (2,4)-- (1,0);
\draw [line width=2pt] (2,4)-- (0,0);
\draw [line width=2pt] (2,4)-- (3,0);
\draw [line width=2pt] (3,4)-- (2,0);
\draw [line width=2pt] (3,4)-- (1,0);
\draw [line width=2pt] (3,4)-- (0,0);
\draw [line width=2pt] (4,4)-- (5,0);
\draw [line width=2pt] (4,4)-- (6,0);
\draw [line width=2pt] (4,4)-- (7,0);
\draw [line width=2pt] (5,4)-- (4,0);
\draw [line width=2pt] (5,4)-- (6,0);
\draw [line width=2pt] (5,4)-- (7,0);
\draw [line width=2pt] (6,4)-- (7,0);
\draw [line width=2pt] (6,4)-- (5,0);
\draw [line width=2pt] (6,4)-- (4,0);
\draw [line width=2pt] (7,4)-- (6,0);
\draw [line width=2pt] (7,4)-- (5,0);
\draw [line width=2pt] (7,4)-- (4,0);
\draw [line width=2pt] (0,4)-- (4,0);
\draw [line width=2pt] (1,4)-- (5,0);
\draw [line width=2pt] (2,4)-- (6,0);
\draw [line width=2pt] (3,4)-- (7,0);
\draw [line width=2pt] (0,0)-- (4,4);
\draw [line width=2pt] (1,0)-- (5,4);
\draw [line width=2pt] (2,0)-- (6,4);
\draw [line width=2pt] (3,0)-- (7,4);
\draw (-0.6,4.7) node[anchor=north west] {$(0,0)$};
\draw (0.4,4.7) node[anchor=north west] {$(0,1)$};
\draw (1.4,4.7) node[anchor=north west] {$(0,2)$};
\draw (2.4,4.7) node[anchor=north west] {$(0,3)$};
\draw (6.4,4.7) node[anchor=north west] {$(1,0)$};
\draw (5.4,4.7) node[anchor=north west] {$(1,1)$};
\draw (4.4,4.7) node[anchor=north west] {$(1,2)$};
\draw (3.4,4.7) node[anchor=north west] {$(1,3)$};
\draw (-0.6,0) node[anchor=north west] {$(2,0)$};
\draw (0.4,0) node[anchor=north west] {$(2,1)$};
\draw (1.4,0) node[anchor=north west] {$(2,2)$};
\draw (2.4,0) node[anchor=north west] {$(2,3)$};
\draw (6.4,0) node[anchor=north west] {$(3,0)$};
\draw (5.4,0) node[anchor=north west] {$(3,1)$};
\draw (4.4,0) node[anchor=north west] {$(3,2)$};
\draw (3.4,0) node[anchor=north west] {$(3,3)$};
\begin{scriptsize}
\draw [fill=black] (0,4) circle (2.5pt);
\draw [fill=black] (1,4) circle (2.5pt);
\draw [fill=black] (2,4) circle (2.5pt);
\draw [fill=black] (3,4) circle (2.5pt);
\draw [fill=black] (4,4) circle (2.5pt);
\draw [fill=black] (5,4) circle (2.5pt);
\draw [fill=black] (6,4) circle (2.5pt);
\draw [fill=black] (7,4) circle (2.5pt);
\draw [fill=black] (0,0) circle (2.5pt);
\draw [fill=black] (1,0) circle (2.5pt);
\draw [fill=black] (2,0) circle (2.5pt);
\draw [fill=black] (3,0) circle (2.5pt);
\draw [fill=black] (4,0) circle (2.5pt);
\draw [fill=black] (5,0) circle (2.5pt);
\draw [fill=black] (6,0) circle (2.5pt);
\draw [fill=black] (7,0) circle (2.5pt);
\end{scriptsize}
\end{tikzpicture}
\caption{Orthogonal Colouring of $H(4,2)$}
\label{Figure: Orthogonal Q4}
\end{figure}

It remains to find orthogonal colourings of Hamming graphs when $q= 2,6$. In the case when $q=2$, $O\chi(H(4,2))=4$, as shown in Figure \ref{Figure: Orthogonal Q4}. Theorem \ref{Theorem: Cartesian Product} and induction gives the following corollary.

\begin{cor}
If $d\geq 1$, then $O\chi(H(4d,2))=2^{2d}$
\end{cor}

There are some open problems left to be explored for orthogonal colourings of Hamming graphs. First, to complete the categorization of even Hamming graphs, an orthogonal colouring of $H(4,6)$ using 36 colours would need to be found. Secondly, orthogonal colourings of odd Hamming graphs, $H(2d+1,q)$, have not been studied. Lastly, multiple orthogonal colourings of Hamming graphs have not been studied. In the case of $H(2,q)$, this would correspond to a collection of mutually orthogonal Latin squares.

\section*{Acknowledgements}

This research was supported by an NSERC Discovery grant.
\bibliographystyle{amsplain}
\bibliography{References}

\end{document}